\documentclass[10pt, reqno]{amsart}
\usepackage{amssymb,latexsym}
\usepackage{eucal}

\newtheorem*{proposition*}{}
\newtheorem*{Lem1}{Lemma 1}
\newtheorem*{Lem2}{Lemma 2}
\newtheorem*{Lem3}{Lemma 3}
\newtheorem*{Prop1}{Proposition 1}
\newtheorem*{Prop2}{Proposition 2}

\newtheorem*{Lem17.3}{Lemma 17.3}
\newtheorem*{Lem18.3}{Lemma 18.3}
\newtheorem*{Lem18.5}{Lemma 18.5}
\newtheorem*{Lem19.0}{Lemma 19.0}
\newtheorem*{Lem19.1}{Lemma 19.1}
\newtheorem*{Lem19.1.3}{Lemma 19.1.3}

\newtheorem*{Thm}{Theorem}
\newtheorem*{Cor}{Corollary}

\newtheorem*{Prob1}{Problem 1}
\newtheorem*{Prob2}{Problem 2}
\newtheorem*{Prob3}{Problem 3}

\DeclareMathOperator{\qv}{\textup{Qvar}}
\DeclareMathOperator{\SQ}{\textup{SQ}}

\newcommand{\e}{{\varepsilon}}

\newcommand{\LL}{\mathcal{L}}
\newcommand{\U}{\mathcal{U}}
\newcommand{\K}{\mathcal{K}}
\newcommand{\G}{\mathcal{G}}

\newcommand{\A}{\mathcal{A}}

\newcommand{\F}{\mathcal{F}}
\newcommand{\X}{\mathcal{X}}
\newcommand{\XL}{\mathcal{X}_{\ell}}

\newcommand{\p}{\partial}
\newcommand{\E}{\mathcal{E}}

\newcommand{\ph}{\varphi}

\begin{document}

\title[On quasivarities of groups and equations over groups]
{On quasivarities of groups and equations over groups}
\author{S.V. Ivanov}
\address{Department of Mathematics\\
University of Illinois \\
Urbana \\   IL 61801 \\ U.S.A.}
\email{ivanov@math.uiuc.edu}
\thanks{Supported in part by NSF grants  DMS 98-01500, DMS 00-99612}
\subjclass[2000]{Primary   20E06, 20F05, 20F06, 20F50}

\begin{abstract}
We prove that the quasivariety of groups generated by finite and locally
indicable groups does not contain the class of periodic groups. This result
is related to (and inspired by) the solvability of equations over groups. The
proof uses the Feit-Thompson theorem on the solvability of finite groups of
odd order, Kostrikin-Zelmanov results on the restricted Burnside problem and
applies techninal details of a recent construction of weakly finitely
presented periodic groups.
\end{abstract}
\maketitle

Let $\X = \{ x_1, x_2, \dots \}$  be a countably
infinite alphabet and  $U_1, \dots,  U_n, V$ be words
in $\X^{\pm 1}= \X \cup \X^{-1}$ (called $\X$-{\em words}).
A {\em quasiidentity} is an expression of the form
$$
( U_1 = 1 \wedge \dots  \wedge U_m = 1 ) \Rightarrow   V = 1 ,
$$
where $\wedge$ and  $\Rightarrow$ are the signs of conjunction and
implication, respectively. A quasiidentity holds in a group $G$ if it is a
true formula for any substitution $g_i \to x_i$, where $g_i \in G$, $i = 1,2
\dots$. A {\em quasivariety} of groups is the class of groups defined by a
set of quasiidentities; that is, the class of all groups in which every
quasiidentity of a given set holds (see \cite{M70} for more details). For
example, the quasiidentity $x^2=1 \to x=1$ defines the class of groups
without involutions.

Recall that an {\em equation} over a group $G$ with an unknown $y$ is an
expression of the form
\begin{equation}
w(y, G) = y^{\e_1} g_1 \dots  y^{\e_\ell} g_\ell = 1 ,
\end{equation}
where $\e_1, \dots, \e_\ell \in \{ \pm 1\}$ and $g_1,  \dots, g_\ell \in G$.
The elements $g_1,  \dots, g_\ell$ of $G$ are called the {\em coefficients}
and $\ell$ is the {\em length} of equation (1). The equation (1)  is called
{\em solvable } if there is a group $H$ that contains $G$ as its subgroup and
there is an element $h \in H$ so that the substitution $h \to y$ in (1)
results in the equality $w(h, G) = 1$ in $H$. Clearly, the equation (1) is
solvable if and only if $G$ naturally embeds in the quotient of the free
product $\langle y \rangle_\infty * G$ by  relation (1), where $\langle y
\rangle_\infty$ is the infinite cyclic group generated by $y$. The equation
(1) is called {\em nonsingular } if $\e_1 + \dots + \e_\ell \neq 0$
(otherwise, it is {\em singular}).

It is a well-known conjecture of Kervaire and Laudenbach that every
nonsingular equation is solvable over any group $G$. We will say that the
equation (1) is {\em strongly solvable } over a group $G$ if it is solvable
for an arbitrary $\ell$-tuple $(g_1, \dots, g_\ell)$ of elements of $G$.
Clearly, the Kervaire-Laudenbach conjecture is equivalent to its strong
version, that is, the Kervaire-Laudenbach conjecture holds if and only if
every nonsingular equation   is strongly solvable over any group $G$. It is
also easy to see that if the  equation (1) is strongly solvable over a
nontrivial group $G$ then this equation is nonsingular.

The following establishes a link between foregoing definitions.

\begin{Prop1} Let $\K(w(y, G))$ denote the class of groups over
which equation $(1)$ is strongly solvable. Then $\K(w(y, G))$ is a
quasivariety of groups whose quasiidentities can be written using $\ell$
letters $x_1, \dots, x_\ell$, where $\ell$ is the length of equation $(1)$.
Furthermore, the class of groups over which all nonsingular equations are
solvable is a quasivariety of groups (denoted by $\qv\E$).
\end{Prop1}

It was proved by Gerstenhaber and  Rothaus \cite{GR62} that the
Kervaire-Laudenbach  conjecture holds for finite groups, and so every finite
group is in $\qv\E$. By Proposition 1, this result can be stated as follows:
$\qv \F \subseteq \qv \E$, where $\qv\F$ is the quasivariety generated by all
finite groups ($\qv\F$ consists of all those groups in which quasiidentities
of finite groups hold).

It was asked by Klyachko whether every periodic group $G$ is in $\qv\F$; this
would imply that every nonsingular equation is solvable over $G$. Note that
any locally finite group, the finitely generated infinite periodic groups
constructed by Golod \cite{Gl64} and the subsequent examples of  Aleshin
\cite{Al72}, Sushchanskii  \cite{S79}, Grigorchuk \cite{Gr80},  Gupta and
Sidki \cite{GS83}, \cite{Gp89}, are all in $\qv\F$, because they are either
locally finite or residually finite.

We also mention that, in view of Proposition 1, one  possible way of proving
the Kervaire-Laudenbach conjecture would be to prove it for groups in a class
$\K$ which generates the quasivariety of all groups. In this regard, recall
that Howie  \cite{H81}  (see also Brodski's articles  \cite{B84},
\cite{B80a}, \cite{B80b}) proved the Kervaire-Laudenbach conjecture for
locally indicable groups (a group $G$ is called {\em locally indicable} if
every nontrivial finitely generated subgroup $H$ of $G$ has an infinite
cyclic homomorphic image). Let $\LL$ denote the class of locally indicable
groups. By Proposition 1, $\qv\LL \subseteq \qv\E$. Since it is also true
that $\qv(\F \cup \LL) \subseteq \qv\E$, it follows that the
Kervaire-Laudenbach conjecture would be proven if one could show that $\qv(\F
\cup \LL)$ is the quasivariety of all groups. To address these natural
questions, we will prove the following.

\begin{Thm}  Let $n >2^{16}$ be odd and $\K$ be a class of
groups such that if $A \in \K$ and $B$
is a finitely generated subgroup of exponent $n$
in $A$ then $B$ is finite (in particular, $\K = \F$ or
$\K = \F \cup \LL$). Then there is an infinite
finitely generated group of exponent $n^3$ which is not
contained in the quasivariety $\qv\K$ of groups generated by $\K$.
\end{Thm}

To prove this Theorem we will use some of technical results of article
\cite{I02} whose construction is based upon \cite{I94} (see Lemmas 1, 2
below). We will also use fundamental results of Feit and Thompson \cite{FT63}
on the solvability of finite groups of odd order and those of Kostrikin
\cite{K59, K86} and Zelmanov \cite{Z91} on the restricted Burnside problem
for finite groups of prime power exponent (together with the Hall-Higman
reduction \cite{HH56} and the classification of finite simple groups; for $n
= p^{\alpha}$, where $p$ is a prime, our proof becomes simpler and uses only
Kostrikin-Zelmanov results  \cite{K59, K86}, \cite{Z91}).

As an immediate corollary of this Theorem, we have negative answers
to both Klyachko's question and the question whether
every group is in $\qv(\F \cup \LL)$.

\begin{Cor}
The quasivariety $\qv(\F\cup\LL)$ generated by finite groups and locally
indicable groups does not contain the class of periodic groups.
\end{Cor}

We remark in passing that no  results are known that deal with the
solvability of nonsingular equations over free Burnside groups of exponent $n
\gg 1$, Tarski monsters and other nonlocally finite groups of bounded
exponent (see \cite{Ad75}, \cite{O89}, \cite{IO91}, \cite{I94}). It is not
even known whether there are such groups in $\qv\E$. In this connection,
observe that every group $G$ has a largest homomorphic image $G/\E(G) \in
\qv\E$, where $\E(G)$ is the normal subgroup of $G$ generated by values of
the  right parts of the quasiidentities of $\qv\E$ whose left parts are true
in $G$. Therefore, the question whether there exists an infinite
$m$-generator group of  exponent $n$ in $\qv\E$ (resp. in $\qv\F$)  is
equivalent to the question whether the quotient $B(m,n)/\E(B(m,n))$ (resp.
$B(m,n)/\F(B(m,n))$) is infinite, where $B(m,n)$ is a free $m$-generator
Burnside group of exponent $n$, that is, $B(m,n) = F_m/F_m^n$ and $F_m$ is a
free group of rank $m$. These problems seem to be rather interesting new
versions of the Burnside problem (which have some similarity with the
restricted Burnside problem, see also Problem 3).

\begin{Prob1}
Does there exist an infinite $m$-generator group of (odd) exponent $n \gg~1$
over which all nonsingular equations are solvable? Equivalently, is the
quotient $B(m,n)/\E(B(m,n))$ infinite?
\end{Prob1}

\begin{Prob2}
Does there exist an infinite $m$-generator group of (odd) exponent $n \gg~1$
in which all quasiidentities of finite groups hold? Equivalently, is the
quotient $B(m,n)/\F(B(m,n))$ infinite?
\end{Prob2}

Since $\qv \F \subseteq \qv \E$, we have that $\E(B(m,n)) \subseteq
\F(B(m,n))$ and so a negative solution to Problem 1 would imply a negative
solution to Problem 2 (we would conjecture, though, that $B(m,n) \in \qv \E$
for (odd) $n \gg 1$).

It is relevant to mention another problem (stated in \cite{KN92},
\cite{IO96}).

\begin{Prob3}
Does there exist an integer $r = r(m,n)>0$, where $m >1$ and $n$ is odd, such
that if $G$ is a finite group, generated by elements $g_1, \dots, g_m$, and
every word $v$ in $g_1^{\pm 1}, \dots, g_m^{\pm 1}$ of length at most $r$
satisfies $v^n =1$ in $G$, then $G$ has exponent $n$?
\end{Prob3}

Recall that a positive solution to Problem 3 for (odd) $n \gg1$ implies the
existence of Gromov hyperbolic groups \cite{G87} that are not residually
finite.

As in the proof of Theorem, putting together results of \cite{FT63} and
\cite{K59, K86}, \cite{Z91}, \cite{HH56} (together with the classification of
finite simple groups), we will show the following.

\begin{Prop2} For given $m$ and $n$, a positive solution to Problem 3 implies
a negative solution to Problem 2.
\end{Prop2}
\medskip

{\em Proof of Theorem.} \ Let $I$ be a finite index set with $|I| > 1$.
Consider an alphabet
$$
\U = \{  a_i, b_i, x_i, y, c \; | \; i \in I \}
$$
and the following relations (2)--(8) in letters of
$\U^{\pm 1} = \U \cup \U^{-1}$.

\begin{gather}
x_i c x_i^{-1} =  c b_i , \  \ i \in I ,   \label{2} \\
y b_i  y^{-1} =   b_i  a_i , \  \ i \in I ,    \label{3} \\
c^n =1 ,  \label{4} \\
x_i b_j  =  b_j  x_i  , \  \ i, j  \in I ,  \label{5} \\
a_i b_j  =   b_j a_i  , \  \ i, j  \in I , \label{6} \\
a_i  c  =   c a_i  , \  \ i   \in I ,  \label{7} \\
y c  =   c y , \  \ i   \in I .  \label{8}
\end{gather}

Let a group $\G$ be given by  the presentation whose alphabet is $\U$ and
whose set of defining relations consists of  relations (2)--(8), thus
\begin{equation} \label{9}
\G = \langle \U \;  \|\  \; (2)-(8)    \rangle .
\end{equation}

Let $F(\A)$ denote a free group in the alphabet $\A$ and $B(\A,n) =
F(\A)/F(\A)^n$  denote the free Burnside group of exponent $n$ in the
alphabet $\A$.

The following two lemmas are proved in \cite{I02}.

\begin{Lem1} The subgroup $\langle   a_i \; | \; i \in I  \rangle$
of the group $\G$,  generated by the images of letters $a_i, i \in I$, is
naturally isomorphic to $B(\A,n)$.
\end{Lem1}

\begin{Lem2}
Let $n > 2^{16}$ be odd. Then the subgroup $\langle   a_i \; | \; i \in I
\rangle$ of the group $\G/\G^{n^3}$, generated by the images of letters
$a_i$, $i\in I$, is naturally isomorphic to $B(\A,n)$. In particular, this
subgroup $\langle  a_i \; | \; i \in I  \rangle$ is infinite.
\end{Lem2}


According to Feit and Thompson \cite{FT63}, every finite group of odd order
is solvable. Furthermore, if $F$ is a finite $m$-generator group of odd
exponent $n$, then Kostrikin-Zelmanov results \cite{K59, K86}, \cite{Z91} on
the restricted Burnside problem, a reduction theorem due to Hall and Higman
\cite{HH56} and the classification of finite simple groups imply that there
is an upper bound $c(m,n)$  for the class of solvability of $F$.

Let us again consider relations (2)--(8) and rewrite each of them in the form
$W = 1$. Construct a quasiidentity whose premise is the conjunction of all
these equalities $W=1$ and whose conclusion is $S_d(v_1, \dots, v_{2^d}) =
1$, where $S_d$ is a "solvable" commutator of class $d$ (so that the quotient
$F(\A)/S_d(F(\A))$ is the free solvable group of class $d$ in the alphabet
$\A$) and $v_1, \dots, v_{2^d}$ are arbitrary elements of $F(\A)$. Picking $d
= c(|I|, n)$, it is easy to see from Lemma 1 that the system of these
quasiidentities (over all possible tuples $(v_1, \dots, v_{2^d})$ of elements
of $F(\A)$) holds in any group $A \in \K$. On the other hand, if this system
holds in the group $\G/\G^{n^3}$, where $\G$ is given by presentation (9),
then the subgroup $\langle   a_i \; | \; i \in I  \rangle$ of $\G/\G^{n^3}$
would have to be solvable of class $d$ and so it would be finite. This
contradiction to Lemma 2 shows that $\G/\G^{n^3}  \not\in \qv\K$ and Theorem
is proved. \qed
\smallskip

Note that if $n = p^{\alpha}$, where $p$ is a prime, then the above argument
becomes simpler, because a finite group of such an exponent $n$ is nilpotent
(so we do not need to refer to  \cite{FT63}) and Kostrikin-Zelmanov results
\cite{K59, K86}, \cite{Z91} directly yield the required upper bound $c(m,n)$
(so we do not need to refer to \cite{HH56} and use the classification  of
finite simple groups).
\smallskip

{\em Proof of Proposition 1.} \ If (1) is a singular equation then $\K(w(y,
G))$ contains the trivial group only and our claim  is obvious. Hence we can
assume that the equation (1) is nonsingular.

Let
$$
\XL = \{ x_1, \dots, x_\ell \} , \quad
w(y, \XL ) = y^{\e_1} x_1 \dots  y^{\e_\ell} x_\ell
$$
and consider a group presentation
\begin{equation} \label{10}
\langle    x_1, \dots, x_\ell, y  \;  \|\  \;
u_1(\XL)=1, \dots, u_k(\XL)=1, \;  w(y, \XL) = 1   \rangle ,
\end{equation}
where $u_1(\XL), \dots, u_k(\XL)$ are words in $\XL^{\pm 1} = \XL \cup \XL^{-1}$,
$k \ge 1$.

A spherical diagram $\Delta$ over presentation  (10) is called {\em $w(y,
\XL)$-regular} if the following properties (W1)--(W3) hold (for the
definition of a (spherical) diagram over a group presentation the reader is
referred to \cite{O89}; see also \cite{LS77}).

\begin{enumerate}
\item[(W1)] If $\Pi$ is a 2-cell of $\Delta$
and $\p \Pi$ is the positively oriented boundary of $\Pi$ then
the word $\ph(\p \Pi)$, up to a cyclic permutation, is either
$w(y, \XL)^{\pm 1}$ (in which case $\Pi$ is called a {\em $y$-face})
or $u_j(\XL)$, where $j \in \{ 1, \dots, k\}$, and for every
$j \in \{ 1, \dots, k\}$ there is precisely
one 2-cell $\Pi$ in $\Delta$ with   $\ph(\p \Pi) =  u_j(\XL)$
(in which case $\Pi$ is called a {\em $u_j$-face}).

\item[(W2)] If $e \in \p \Pi$ is an (oriented) edge and $\Pi$ is
a $u_j$-face in $\Delta$ then $e^{-1} \in \p \Pi'$, where $\Pi'$
is a $y$-face.

\item[(W3)] Every vertex $o$ in $\Delta$ has degree 3 and if $e_1, e_2, e_3$
are 3 (oriented) edges whose terminal vertex is $o$
then for one of them, say, $e_1$, one has $\ph(e_1) = y^{\pm 1}$ and
$e_2, e_3^{-1}$ are consecutive edges of $\p \Pi^{\pm 1}$, where $\Pi$ is
a $u_j$-face.
\end{enumerate}

For every $w(y, \XL)$-regular spherical diagram $\Delta$ we consider
$k$ quasiidentities
$$
\left(  u_{1}(\XL) =1  \wedge  \dots \wedge
u_{j-1}(\XL) =1 \wedge u_{j+1}(\XL) =1   \wedge \dots
\wedge u_{k}(\XL) =1  \right)  \Rightarrow  u_{j}(\XL) =1
$$
for $j =1, \dots, k$.

By $\SQ(w(y, \XL))$ denote the system of such  quasiidentities over all
$w(y, \XL)$-regular spherical diagrams $\Delta$ (for all $k =1,2,\dots$).
Finally, let  $\qv(\SQ(w(y, \XL)))$ denote the quasivariety defined by the system
$\SQ(w(y, \XL))$.

\begin{Lem3} The class $\K(w(y, G))$
of groups over which  the equation $(1)$ is strongly solvable is
$\qv(\SQ(w(y, \XL)))$.
\end{Lem3}

{\em Proof.} To prove that $\K(w(y, G)) \subseteq  \qv(\SQ(w(y, \XL)))$,
suppose, on the contrary,  that $G_1$ is a group such that
$G_1 \in \K(w(y, G))$ and $G_1 \not\in \qv(\SQ(w(y, \XL)))$.

Since $G_1 \not\in \qv(\SQ(w(y, \XL)))$, there is a quasiidentity
$$
\left(  u_{1}(\XL) =1  \wedge  \dots \wedge u_{j-1}(\XL) =1
\wedge u_{j+1}(\XL) =1  \wedge \dots \wedge u_{k}(\XL) =1  \right)
\Rightarrow  u_{j}(\XL) =1
$$
in the system  $\SQ(w(y, \XL))$ which does not hold in $G_1$.

Assume that this quasiidentity fails for an $\ell$-tuple
$(g_{1,1}, \dots, g_{1, \ell})$ of elements of $G_1$ and
let $\Delta$ be the $w(y, \XL)$-regular spherical diagram
used in construction of this quasiidentity.
Relabelling edges of $\Delta$ whose labels are
$x_1^{\pm 1}, \dots, x^{\pm 1}_\ell$ by corresponding elements
$g_{1,1}^{\pm 1}, \dots, g^{\pm 1}_{1,\ell}$ and taking
$u_j$-face out of $\Delta$, we will get a disk diagram
$\Delta_0$ over the relative presentation
$$
\langle   y, G_1  \  \|\  \
w(y, g_{1,1}, \dots, g_{1,\ell} ) = 1  \rangle
$$
so that
$\ph(\p \Delta_0) = u_j( g_{1,1}, \dots, g_{1,\ell}) \neq 1$
in $G_1$. This proves that equation (1) is not strongly solvable over
$G_1$, contrary to $G_1 \in \K(w(y, G))$.

To prove the converse inclusion  $\qv(\SQ(w(y, \XL))) \subseteq  \K(w(y, G)) $,
assume, on the contrary, that $G_1 \in \qv(\SQ(w(y, \XL)))$ and
$G_1 \not\in \K(w(y, G))$.

Since $G_1 \not\in \K(w(y, G))$, there is an $\ell$-tuple
$(g_{1,1}, \dots, g_{1, \ell})$ of elements of $G_1$ such that
the equation
$
w(y, g_{1,1}, \dots, g_{1, \ell} )
= y^{\e_1} g_{1,1} \dots  y^{\e_\ell} g_{1,\ell} = 1
$
is not solvable. This means the existence of a disk diagram $\Delta_0$
over the relative presentation
$
\langle  y, G_1  \  \|\  \
w(y, g_{1,1}, \dots, g_{1,\ell} ) = 1  \rangle
$
so that $\ph(\p \Delta_0) \in G_1$ and  $\ph(\p \Delta_0) \neq 1$ in $G_1$.
Now it is easy to construct a $w(y, \XL)$-regular
spherical diagram $\Delta$ from $\Delta_0$ such that one of
quasiidentities associated as above with $\Delta$ does not hold
for the $\ell$-tuple $(g_{1,1}, \dots, g_{1, \ell})$. This contradiction
to $G_1 \in \qv(\SQ(w(y, \XL)))$ completes the proof of Lemma 3.  \qed
\smallskip

The first claim of Proposition 1 follows from definitions and Lemma 3.
Furthermore, in view of Lemma 3, the union $\bigcup \SQ(w(y, \XL))$ over all
nonsingular equations  (1) is a system of quasiidentities  which defines the
class of groups over which  all nonsingular equations are strongly solvable.
Therefore, this class is a quasivariety of groups and Proposition 1 is
proved. \qed
\smallskip

{\em Proof of Proposition 2.} \ Assume that Problem 3 is solved in the
affirmative for given $m, n$ and $r(m,n)$ is a required integer. As in the
proof of Theorem, it follows from results of \cite{FT63}, \cite{K59, K86},
\cite{Z91}, \cite{HH56} (and the classification of finite simple groups) that
there is a bound $c(m,n)$ for the class of solvability  of an $m$-generator
finite group of odd exponent $n$. As in the proof of Theorem, consider a
"solvable" commutator $S_d$ of class $d = c(m,n)$. Define a system of
quasiidentities whose premise is the conjunction of equalities $u(x_1, \dots,
x_m)^n =1$, where the word $u(x_1, \dots, x_m)$ has length at most $r(m,n)$,
and whose conclusion is $S_d(v_1, \dots, v_{2^d}) = 1$, where $v_1, \dots,
v_{2^d}$ are arbitrary words in $\{ x_1^{\pm 1}, \dots,  x_m^{\pm 1} \}$. It
is easy to see that every quasiidentity of this system  holds in any finite
group. Therefore, the subgroup $\F(B(m,n))$  contains the verbal subgroup
$S_d(B(m,n))$ and so the quotient $B(m,n)/ \F(B(m,n))$  is solvable and
finite. This proves Proposition~2. \qed

\end{document}